\documentclass{amsart}
\usepackage[utf8]{inputenc}
\usepackage{amsmath}
\usepackage{hyperref}
\usepackage{dynkin-diagrams}
\usepackage[shortlabels]{enumitem}

\newtheorem{theorem}{Theorem}[section]
\newtheorem{lemma}[theorem]{Lemma}

\theoremstyle{definition}
\newtheorem{definition}[theorem]{Definition}
\newtheorem{example}[theorem]{Example}

\theoremstyle{remark}
\newtheorem{remark}[theorem]{Remark}

\numberwithin{equation}{section}

\DeclareMathOperator{\GH}{Gh}
\DeclareMathOperator{\Tr}{Tr}
\DeclareMathOperator{\ann}{ann}
\DeclareMathOperator{\PL}{PL}
\DeclareMathOperator{\Herm}{Herm}
\DeclareMathOperator{\Rerm}{Re}
\DeclareMathOperator{\PSU}{PSU}
\DeclareMathOperator{\PSL}{PSL}
\DeclareMathOperator{\Co}{Co}
\DeclareMathOperator{\Suz}{Suz}
\DeclareMathOperator{\HJ}{HJ}

\begin{document}

\author{Benjamin Nasmith}

\address{Department of Mathematics and Computer Science \\
Royal Military College of Canada \\
Kingston, ON}

\email{Ben.Nasmith@gmail.com}

\thanks{\textit{Acknowledgments}. The published version of this pre-print is available as \textit{Algebraic Combinatorics}, Volume 5 (2022) no. 3, pp. 401-411. doi : 10.5802/alco.215. \href{https://alco.centre-mersenne.org/articles/10.5802/alco.215/}{https://alco.centre-mersenne.org/articles/10.5802/alco.215/} }

%%%%% Sujet

\keywords{Projective t-designs, octonions, Leech lattice, generalized hexagon}
 
\subjclass{05B99, 17A75, 17C40}

%%%%% Gestion

% \DOI{10.5802/alco.215}
% \datereceived{2020-09-03}
% \daterevised{2021-12-15}
% \dateaccepted{2021-12-19}

%%%%% Titre et rÃ©sumÃ©
\title{Octonions and the two strictly projective tight 5-designs}

\begin{abstract}
 In addition to the vertices of the regular hexagon and icosahedron, there are precisely two strictly projective tight 5-designs: one constructed from the short vectors of the Leech lattice and the other corresponding to a generalized hexagon structure in the octonion projective plane. 
 This paper describes a new connection between these two strictly projective tight 5-designs---a common construction using octonions.
 Certain octonion involutionary matrices act on a three-dimensional octonion vector space to produce the first 5-design and these same matrices act on the octonion projective plane to produce the second 5-design. 
 This result uses the octonion construction of the Leech lattice due to Robert Wilson and provides a new link between the generalized hexagon Gh(2,8) and the Leech lattice.
\end{abstract}

%%%%%%%%%%%%%%%%%%%%%%%%%%%%%%%%%%%%%%%%%%%

\maketitle

\section{Introduction}

A spherical $t$-design is a finite subset of points on the unit sphere in a real vector space with the following special property: the average value of any polynomial of degree at most $t$ over the sphere is equal to the average value of the polynomial evaluated at the points of the $t$-design~\cite{delsarte_spherical_1977}. 
Every $t$-design is also an $A$-code, where $A$ is the set of angles between distinct points in the $t$-design.
A projective $t$-design (and $A$-code) generalizes this concept from spheres to projective spaces~\cite{neumaier_combinatorial_1981, hoggar_t-designs_1982}.
Taken together, the spheres and infinite projective spaces constitute the compact symmetric spaces of rank $1$~\cite{hoggar_t-designs_1992}, the compact and connected two-point homogeneous spaces~\cite{wang_two-point_1952}, and also the manifolds of primitive idempotents for the simple Euclidean Jordan algebras~\cite{faraut_analysis_1994}. 
That is, for $V$ a simple Euclidean Jordan algebra of rank $\rho$ and degree $d$, the manifold of primitive idempotents $\mathcal{J}(V)$ is given by,
\begin{align*}
 \Omega_{d+1}, \quad 
 \mathbb{RP}^{\rho - 1}, \quad 
 \mathbb{CP}^{\rho-1}, \quad 
 \mathbb{HP}^{\rho-1}, \quad 
 \mathbb{OP}^{2}, \quad d \ge 1,~\rho \ge 3.
\end{align*}
There are no repetitions on this list, but the projective lines (not listed) are isomorphic to the following spheres:
\begin{align*}
 \Omega_2 \cong \mathbb{RP}^1, \quad 
 \Omega_3 \cong \mathbb{CP}^1, \quad 
 \Omega_5 \cong \mathbb{HP}^1, \quad 
 \Omega_9 \cong \mathbb{OP}^1.
\end{align*}
A $t$-design is \textit{spherical} when it is a subset of sphere $\Omega_{d+1}$ (with $\rho = 2$), and \textit{projective} when it is a subset of projective space $\mathbb{FP}^{\rho - 1}$ (with $d = [\mathbb{F}:\mathbb{R}] = 1,2,4,8$). We will call a $t$-design \textit{strictly projective} when it has $\rho \ge 3$ and is therefore not also spherical. 
Both spherical and projective $t$-designs are interesting objects in part because of their connections to real, complex, and quaternion reflection groups, as well as other sporadic simple groups~\cite{shephard_finite_1954, cohen_finite_1976, cohen_finite_1980, hoggar_t-designs_1982, bannai_survey_2009}.

A \emph{tight} $t$-design (whether spherical or projective) simultaneously meets a \emph{lower bound}, given its value of $t$, and an \emph{upper bound}, given its value of $A$.
While $t$-designs are common, tight $t$-designs are rare and continue to elude full classification. 
Even so, a projective tight $t$-design with $\mathbb{FP}^{\rho-1} \ne \mathbb{RP}^1$ has $t \le 5$~\cite{hoggar_tight_1984, hoggar_tight_1989} and there are precisely four projective tight $5$-designs. Two of these are also spherical, the hexagon in $\Omega_2 \cong \mathbb{RP}^1$ and the icosahedron in $\Omega_3 \cong \mathbb{CP}^1$~\cite{lyubich_tight_2009}. 
The remaining two tight $5$-designs are strictly projective: one in $\mathbb{RP}^{23}$ consisting of the $98280$ lines spanned by the Leech lattice short vectors and the other in $\mathbb{OP}^2$ realizing the $819$ lines of the unique generalized hexagon $\GH(2,8)$~\cite{hoggar_tight_1989}. 

Hoggar conjectured that the two strictly projective tight $5$-designs are closely related, since the first has cardinality $98280 = 120\cdot 819$, the second has cardinality $819$, and there are $120$ pairs of opposite octonion integer units. Hoggar reports that this conjecture was initially met with skepticism~\cite{hoggar_t-designs_1982}.
This paper provides the missing common construction using certain octonion involutionary matrices. 
These matrices can act on both the octonion vector space $\mathbb{O}^3 \cong \mathbb{R}^{24}$ and the octonion projective plane $\mathbb{OP}^2$ to produce the two strictly projective tight $5$-designs.
This common construction also serves as a new connection between the generalized hexagon $\GH (2,8)$ and the Leech lattice.

\section{Tight t-Designs}

This section reviews tight $t$-designs and their partial classification in order to provide context for the common construction that follows.

\subsection{Definitions}
Let $V$ be a simple Euclidean Jordan algebra of rank $\rho$ and degree $d$ and let $\mathcal{J}(V)$ be the manifold of primitive idempotents. 
We can use the following \textit{renormalized Jacobi polynomials} to describe both the spherical ($\rho = 2$) and projective ($d = 1,2,4,8$) cases:
\begin{align*}
 Q_k^\varepsilon(x) = \left(\frac{1}{2}\rho d + 2k + \varepsilon - 1\right)\frac{(\frac{1}{2}\rho d )_{k +\varepsilon - 1} }{(\frac{1}{2}d)_{k+\varepsilon} }P_k^{(\frac{1}{2}d (\rho - 1) - 1, \frac{1}{2}d - 1 + \varepsilon)}(2x-1).
\end{align*}
Here we have $k$ a non-negative integer, $\varepsilon = 0$ or $1$, Pochhammer symbol $(a)_i = a (a+1)(a+2)\cdots(a+i-1)$, and Jacobi polynomial $P_n^{(\alpha,\beta)}(x)$ as defined in~\cite[22.2.1]{abramowitz_handbook_1972}.
Let $X$ be a finite subset of $\mathcal{J}(V)$ and let $\langle x,y\rangle = \Tr (x\circ y)$ be the Jordan inner product. 
We can compute \textit{angle set} $A$ as follows:
\begin{align*}
 A(X) = \left\{ \langle x,y\rangle \mid x, y \in X \subset \mathcal{J}(V), x \ne y \right\}. 
\end{align*}
We can also compute \textit{strength} $t$ as the maximum non-negative integer $t$ that satisfies%,
\begin{align*}
 \sum_{x \in X} \sum_{y\in X} Q_k^{0}(\langle x,y\rangle) = 0, \quad k = 1,2, \ldots, t.
\end{align*}
With these values we call $X$ both an \textit{$A$-code} and \textit{$t$-design}.

\begin{remark}
 In the spherical $\rho = 2$ case, the Jordan inner product has the form $\langle x,y\rangle = \frac{1}{2} + \frac{1}{2}\cos\theta$, where $\theta$ is the angle between the two points $x,y$ on the sphere. This ensures that antipodal points on the sphere are orthogonal relative to the Jordan inner product. 
 In general, for primitive idempotents $x \ne y$, we have $0 \le \langle x,y\rangle < 1$ and $\langle x,x\rangle = 1$.
\end{remark}

The \textit{annihilator polynomial} of $X \subset \mathcal{J}(V)$, denoted $\ann(x)$, is the unique degree $|A|$ polynomial constructed to ensure that%,
\begin{align*}
 |X| = \ann (1), \quad A(X) = \{ \alpha \in \mathbb{R} \mid \ann (\alpha) = 0\}.
\end{align*}
Given finite $X$, we can compute $t$, $A$, and $\ann (x)$ directly. The more difficult task is to specify $t$, $A$, or $\ann (x)$ and then find a finite subset $X \subset \mathcal{J}(V)$ that realizes them. 
In general a $t$-design (and $A$-code) satisfies the inequality $t \le 2 s - \varepsilon$, where $s = |A(X)|$ and $\varepsilon = |A(X) \cap \{0\}|$. 
A tight $t$-design achieves $t = 2s - \varepsilon$. 
It also simultaneously achieves the \emph{lowest} possible cardinality $|X|$ for the given $t$ value and the \emph{highest} possible cardinality $|X|$ for the given $s = |A|$ value~\cite{hoggar_t-designs_1982, bannai_tight_1985}. 
We can equivalently define a \textit{tight $t$-design} as a finite subset $X \subset \mathcal{J}(V)$ with the annihilator polynomial $\ann (x) = x^\varepsilon R_{s-\varepsilon}^\varepsilon(x)$, where we have%,
\begin{align*}
 R_{s-\varepsilon}^\varepsilon(x) = Q_0^\varepsilon(x) + Q_1^\varepsilon(x) + \cdots + Q_{s-\varepsilon}^\varepsilon(x). 
\end{align*}
That ensures that a tight $(2s - \varepsilon)$-design has cardinality $|X| = R_{s-\varepsilon}^\varepsilon(1)$ and an angle set $A(X)$ given by the roots of polynomial $x^\varepsilon R_{s-\varepsilon}^\varepsilon(x)$.

\subsection{Partial Classification}
The partial classification of tight $t$-designs is as follows. 
In the case of the circle $\Omega_2 \cong \mathbb{FP}^1$ ($\rho = 2$, $d = 1$), a tight $t$-design exists for all positive integer $t$ values, corresponding to the vertices of a $(t+1)$-gon. 

As described in~\cite{colbourn_handbook_2007, bannai_survey_2009}, in the remaining spherical cases ($\rho = 2$, $d \ge 2$) a tight $t$-design must have $t = 1,2,3,4,5,7,11$. 
On the sphere $\Omega_{d+1}$, a tight $1$-design is a pair of antipodal points, a tight $2$-design is a simplex of $d+2$ points, and a tight $3$-design is a cross polytope of $2d + 2$ points. 
Tight spherical $4$- and $5$-designs are in one-to-one correspondence and the search for tight spherical $5$-designs is still open. The only known examples outside of $\Omega_2$ are sets of vectors spanning equiangular lines in $\Omega_3$, $\Omega_{7}$, and $\Omega_{23}$. If any further example exists, it will have $d \ge 118$~\cite{bannai_nonexistence_2005, bannai_survey_2009}.
Likewise, the search for spherical tight $7$-designs is open with examples known in $\Omega_8$ and $\Omega_{23}$. If any further spherical tight $7$-designs exist, they will have $d \ge 103$~\cite{bannai_nonexistence_2005, bannai_survey_2009}.
Finally, there is precisely one spherical tight $11$-design, the points defined by the short vectors of the Leech lattice in $\Omega_{24}$.

As described in~\cite{hoggar_tight_1984, hoggar_tight_1989, lyubich_tight_2009}, the remaining projective cases ($\rho \ge 3$, $d = 1,2,4,8$) have $t = 1,2,3,5$. 
A projective tight $1$-design is a Jordan frame, namely $\rho$ orthogonal points. A Jordan frame exists in each projective space and generalizes antipodal points on the sphere for $\rho = 2$. 
All remaining examples have $t \ge 2$.
A real projective tight $t$-design ($d = 1$, $\rho \ge 3$), corresponds to a spherical tight $(2t+1)$-design, so the search for spherical tight $5$- and $7$-designs is equivalent to the search for real projective tight $2$- and $3$-designs.
There is just one real projective tight $5$-design: the lines spanned by the short vectors of the Leech lattice.
The remaining complex or quaternion projective tight $t$-designs ($d = 2,4$, $\rho \ge 3$) must have $t = 2,3$. The search for complex and quaternion projective tight $2$- and $3$-designs remains open. As described in~\cite{cohn_optimal_2016}, there are many known examples of tight $2$-designs in $\mathbb{CP}^{\rho - 1}$, but surprisingly few have been found in $\mathbb{HP}^{\rho-1}$.
Finally, the remaining tight $t$-designs in the octonion case ($d = 8$, $\rho = 3$) must have $t = 2,5$. 
The tight $2$-design in $\mathbb{OP}^2$ has been proven to exist without an explicit construction in~\cite{cohn_optimal_2016}. 
The tight $5$-design in $\mathbb{OP}^2$ was constructed in~\cite{cohen_exceptional_1983}.

To summarize, the classification of tight $t$-designs will remain open until the tight $2$- and $3$-designs in $\mathbb{RP}^{\rho -1}$, $\mathbb{CP}^{\rho-1}$, and $\mathbb{HP}^{\rho -1}$ are all identified. In contrast, the projective tight $5$-designs are fully classified.

\begin{theorem}
 \label{5designclassification}
 A projective tight $5$-design $X \subset \mathcal{J}(V)$ is either,
 \begin{enumerate}
  \item The vertices of a regular hexagon in $\Omega_2 \cong \mathbb{RP}^1$; 
  \item The vertices of a regular icosahedron in $\Omega_3 \cong \mathbb{CP}^1$;
  \item The lines spanned by the short vectors of the Leech lattice in $\mathbb{RP}^{23}$; or
  \item The unique realization of $\GH (2,8)$ in $\mathbb{OP}^2$.
 \end{enumerate}
\end{theorem}

\begin{remark}
 Two of the four projective tight $5$-designs in Theorem~\ref{5designclassification} are also spherical tight $5$-designs constructed from systems of equiangular lines in $\mathbb{R}^2$ and $\mathbb{R}^3$. The other two examples, in $\mathbb{RP}^{23}$ and $\mathbb{OP}^2$, are not spherical and constitute the only two \textit{strictly projective tight $5$-designs}. 
 This paper will identify a common construction for the two unique strictly projective $5$-designs of Theorem~\ref{5designclassification}.
\end{remark}

\begin{remark}
 The proof in~\cite{hoggar_tight_1984, hoggar_tight_1989} that $t = 1,2,3,5$ for projective tight $t$-designs ($d = 1,2,4,8$, $\rho \ge 2$) other than $\mathbb{RP}^1$ ($d = 1$, $\rho = 2$) rests on a faulty lemma.
 Specifically, Hoggar attempts to prove that for a projective tight $t$-design $X$, the angle set $A(X)$ must be rational~\cite{hoggar_tight_1984}, and the proofs of various restrictions on $t$ depend on this result. 
 However, as described in~\cite{lyubich_tight_2009}, the icosahedron vertices in $\mathbb{CP}^1 \cong \Omega_3$ serve as a counter-example since the angle set of that projective tight $5$-design is irrational. 
 Lyubich repairs the faulty lemma in~\cite{hoggar_tight_1984} for $d = 1,2,4$, accounting for the exceptions in $\mathbb{RP}^1 \cong \Omega_2$ and $\mathbb{CP}^1 \cong \Omega_3$, but ignores the octonion $d = 8$ cases~\cite{lyubich_tight_2009}.
 The repair in~\cite{lyubich_tight_2009} involves correctly identifying the idempotent basis of the Bose--Mesner algebra of a tight $t$-design, which was incorrectly chosen in~\cite{hoggar_tight_1984}.
 This can also be done for the remaining spherical ($\rho = 2$) and octonion ($d= 8$) cases in the same way as outlined in~\cite{lyubich_tight_2009}. 
 No new exceptions exist beyond those captured in~\cite{lyubich_tight_2009}.
 As a result, the results about possible $t$ values for tight $t$-designs in~\cite{hoggar_tight_1984, hoggar_tight_1989} still hold true. Furthermore, aside from the known exceptions in $\mathbb{RP}^1 \cong \Omega_2$ and $\mathbb{CP}^1 \cong \Omega_3$, the angle set of a tight $t$-design is indeed rational. A subsequent paper will provide this general proof.
\end{remark}

\section{Octonions and Isometries}

This section describes how certain involutionary isometries of vector space $\mathbb{F}^\rho$ and projective space $\mathbb{FP}^{\rho -1}$, with associative $\mathbb{F} = \mathbb{R}, \mathbb{C}, \mathbb{H}$, generalize to the non-associative octonion case where $\mathbb{F} = \mathbb{O}$.

\subsection{Definitions}
The division composition algebras over the real numbers are precisely the real numbers $\mathbb{R}$, the complex numbers $\mathbb{C}$, the quaternions $\mathbb{H}$, and the octonions $\mathbb{O}$. 
A standard basis for the octonions is $\{i_t \mid t \in \PL (7) = \{\infty\} \cup \mathbb{F}_7\}$, with $1 = i_\infty$ the identity and
\begin{align*}
 i_t^2 = -1, \quad 
 i_t = i_{t+1} i_{t+3} = - i_{t+3} i_{t+1}, \quad 
 t \in \mathbb{F}_7.
\end{align*}
Octonion conjugation is the $\mathbb{R}$-linear involution defined by $\overline{1} = 1$ and $\overline{i}_{t} = -i_t$ for $t \in \mathbb{F}_7$. 
The real-valued norm is given by $N(x) = x\overline{x} = \overline{x} x$. 
The subalgebra of $\mathbb{O}$ generated by a single octonion is commutative (isomorphic to $\mathbb{R}$ or $\mathbb{C}$) and the subalgebra generated by any two octonions is associative (isomorphic to $\mathbb{R}$, $\mathbb{C}$, or $\mathbb{H}$).
Many further details about this non-associative algebra are available in~\cite{springer_octonions_2000, baez_octonions_2002, conway_quaternions_2003, schafer_introduction_2017}.

As described above, simple Euclidean Jordan algebras of rank $\rho \ge 3$ can be described as hermitian matrices relative to octonion conjugation, which we denote $\Herm (\rho,\mathbb{F})$, with the \textit{Jordan product} defined in terms of the usual matrix product $xy$ as follows:
\begin{align*}
 x \circ y = \frac{1}{2}(xy + yx).
\end{align*}
Here we have $\mathbb{F} = \mathbb{R}, \mathbb{C}, \mathbb{H}, \mathbb{O}$ (in the octonion case, we must have $\rho = 3$ and the underlying matrix product $xy$ is non-associative). 
In addition to the Jordan product, each Jordan algebra element defines an endomorphism $P(x)$, known as the \textit{quadratic representation}~\cite[II.3]{faraut_analysis_1994}:
\begin{align*}
 P(x): y \mapsto 2 x\circ(x\circ y) - x^2 \circ y. 
\end{align*}
When the Jordan product is constructed from an \emph{underlying associative algebra} (as in all but the $\rho =3$ and $d = 8$ octonion case) then $P(x)$ simplifies to $P(x)y = xyx$, with $xyx$ computed using that underlying associative product. 
Many further details about Euclidean Jordan algebras are available in~\cite{faraut_analysis_1994, springer_octonions_2000, baez_octonions_2002, schafer_introduction_2017}. 

Let $x = (x_1, x_2, \ldots, x_\rho)$ be a row vector in $\mathbb{F}^\rho$ and let $x^\dagger$ be the conjugate transpose, a column vector. 
We will call a vector $x$ in $\mathbb{F}^\rho$ \textit{commutative} when the coefficients $\{x_1, x_2, \ldots, x_\rho\}$ generate a real or complex subalgebra of $\mathbb{O}$ and \textit{associative} when the coefficients generate a real, complex, or quaternion subalgebra of $\mathbb{O}$. 
We can extend the norm on $\mathbb{F}$ described above to vectors in $\mathbb{F}^\rho$ as follows:
\begin{align*}
 N(x) = x x^\dagger = x_1 \overline{x}_1 + \cdots + x_\rho \overline{x}_{\rho}
 = N(x_1) + \cdots + N(x_\rho).
\end{align*}
This norm is real-valued and serves as a Euclidean norm for the corresponding real vector space $\mathbb{R}^{\rho d} \cong \mathbb{F}^{\rho}$ with $d = [ \mathbb{F}: \mathbb{R}]$. The inner product is constructed from the norm in the standard way: 
\begin{align*}
 (x,y) = \frac{1}{2}\left(N(x+y) - N(x) - N(y)\right) = \frac{1}{2}(x y^\dagger + y x^\dagger) = \Rerm (x y^\dagger) = \Rerm (y x^\dagger).
\end{align*}
If vector $x$ is associative, then $[x] = x^\dagger x / N(x)$ is also a primitive idempotent in $\mathbb{FP}^{\rho-1} \subset \Herm (\rho, \mathbb{F})$. 
Indeed, any primitive idempotent $[x]$ in $\mathbb{FP}^{\rho-1}$ can be constructed this way for some (non-unique) associative vector $x$ in $\mathbb{F}^\rho$.

\subsection{Isometries}

Given our real-valued inner products defined on both $\mathbb{F}^\rho$ and $\mathbb{FP}^{\rho-1}$, we now want to construct isometries.
An important property of the quadratic representation is that when $w\circ w = I_\rho$ the map $P(w)$ is an involutionary automorphism of the Jordan algebra and an isometry of the manifold of primitive idempotents $\mathcal{J}(V)$ relative to $\langle x,y\rangle = \Tr (x\circ y)$.

For associative $\mathbb{F} = \mathbb{R}, \mathbb{C}, \mathbb{H}$, the following pairs of maps defined by associative vector $r$ are involutionary isometries of $\mathbb{F}^\rho$ and $\mathbb{FP}^{\rho-1}$ respectively:
\begin{align}
 \label{pairofisometries}
 x \mapsto x (I_\rho - 2 [r]), && 
 [x] \mapsto (I_\rho - 2 [r]) [x] (I_\rho - 2 [r]).
\end{align}
Matrices of the form $W(r) = I_\rho - 2 [r]$ satisfy $W(r)^\dagger W(r) = I_\rho$ and therefore belong to the matrix groups $O(\rho)$, $U(\rho)$, or $Sp(\rho)$ respectively when $\mathbb{F} = \mathbb{R}, \mathbb{C}, \mathbb{H}$.
In the non-associative case, with $\mathbb{F} = \mathbb{O}$ and $\rho = 3$, we can ensure that the maps of Eq.~\eqref{pairofisometries} are isometries by selecting a \emph{commutative} vector $r$ in $\mathbb{O}^3$. 

\begin{lemma}
 \label{isometrylemma}
 Let $r$ be a commutative vector in $\mathbb{F}^\rho$. Then the maps of Eq.~\eqref{pairofisometries} are respectively isometries of $\mathbb{F}^\rho \cong \mathbb{R}^{\rho d}$ with inner product $(x,y) = \Rerm (x y^\dagger)$ and of $\mathbb{FP}^{\rho-1}$ with inner product $\langle x,y\rangle = \Tr (x\circ y)$. 
\end{lemma}

\begin{proof}
 Consider the map $x \mapsto x W(r)$ with $W(r) = I_\rho - 2 [r]$. In the associative cases with $\mathbb{F} = \mathbb{R}, \mathbb{C}, \mathbb{H}$, the matrix $W(r)$ belongs to one of the matrix groups $O(\rho)$, $U(\rho)$, or $Sp(\rho)$ since $W(r)^2 = I_\rho$. 
 This ensures that the map above is an isometry for $\mathbb{F}$ associative.
 It remains to check the non-associative case with $\mathbb{F} = \mathbb{O}$ and $\rho = 3$, which includes $\rho = 2$ when we restrict to the appropriate subspace.
 Any linear transformation acting on $\mathbb{O}^3$ preserving the norm $N(x)$ will also preserve the inner product $(x,y)$.
 Since our map is linear we need to show that $N(x) = N(x W(r))$.
 To do so, let $x = a + b + c$ with $a = (A,0,0)$, $b = (0,B,0)$, $c = (0,0,C)$ and with $A,B,C$ in $\mathbb{O}$. Writing $W = W(r)$, we have%,
 \begin{align*}
  N(xW) 
  &= N(aW) + N(bW) + N(cW) + 2(a W, bW) %\\
  + 2(b W, cW)+ 2(c W, aW).
 \end{align*}
\looseness-1
 The first term satisfies $N(aW) = N(a)$ since the coefficients of $a$ and $W$ belong to a common quaternion subalgebra. Likewise we have $N(bW) = N(b)$ and $N(cW) = N(c)$. Since $N(x) = N(a) + N(b) + N(c)$, it remains to show that the cross terms of the form $(aW, bW)$ vanish.
 To do so we write primitive idempotent matrix $[r]$ in the form%,
 \begin{align*}
  [r] = \left(\begin{matrix}
   d & F & \overline{E} \\
   \overline{F} & e & D \\
   E & \overline{D} & f
  \end{matrix}\right), ~d,e,f \in \mathbb{R},~D,E,F \in \mathbb{O}.
 \end{align*}
 The following inner product evaluates to%,
 \begin{align*}
  \frac{1}{4}(a W, b W) &= \Rerm ((A \overline{E})(\overline{D}~\overline{B})) +
  \left(e-\frac{1}{2} \right) \Rerm ((AF)\overline{B})
  +\left(d-\frac{1}{2} \right) \Rerm (A(F\overline{B})).
 \end{align*}
 In general, $\Rerm (A(F\overline{B})) = \Rerm ((AF)\overline{B})$ for any octonions $A, B, F$~\cite[145]{wilson_finite_2009}. Likewise, $\Rerm ((A \overline{E})(\overline{D}~\overline{B})) = \Rerm (((A \overline{E})\overline{D})\overline{B})$. We can also use the primitive idempotent relations $e + d - 1 = -f$ and $f F = \overline{E}~\overline{D}$~\cite[157]{wilson_finite_2009}.
 Finally, by construction $E$ and $D$ belong to a common complex subalgebra of $\mathbb{O}$, since $r$ is a commutative vector.
 This ensures that $(A\overline{E})\overline{D} = A(\overline{E}~\overline{D})$. Taken together, our expression simplifies to zero:
 \begin{align*}
  \frac{1}{4}(a W, b W) &= \Rerm (((A \overline{E})\overline{D})\overline{B}) +
  \left(d+e-1 \right) \Rerm ((AF)\overline{B}) \\
  &= \Rerm (((A \overline{E})\overline{D})\overline{B}) - \Rerm ((A (fF))\overline{B}) \\
  &= \Rerm (((A \overline{E})\overline{D} - A(\overline{E}~\overline{D}))\overline{B}) \\
  &= 0.
 \end{align*}
 A similar calculation cycling $a\mapsto b \mapsto c \mapsto a$, $d\mapsto e \mapsto f \mapsto d$, and $D\mapsto E \mapsto F \mapsto D$ verifies that $(bW,cW) = (cW, aW) = 0$. 
 This confirms that $x \mapsto x W(r)$ is an isometry of $\mathbb{O}^3$ when $r$ is a commutative vector.

 In the associative cases, with $\mathbb{F} = \mathbb{R}, \mathbb{C}, \mathbb{H}$, the map $[x] \mapsto W(r)[x]W(r) = P(W(r))[x]$ is a known isometry of $\mathbb{FP}^{\rho-1}$. 
 In the non-associative case with $\mathbb{F} = \mathbb{O}$,
 we verify that $[x] \mapsto W(r)[x]W(r)$ is an isometry of $\mathbb{OP}^2$ by beginning with the known isometry $P(W(r))[x]$ given by the quadratic map. 
 We can write $[x] = a + b + c + A + B + C$ for $a,b,c$ real-valued matrices corresponding to the diagonal entries of $[x]$ and $A,B,C$ octonion-valued Hermitian matrices corresponding to the octonion valued off-diagonal entries of $[x]$. 
 \begin{align*}
  P(W(r))[x] = P(W(r))a + \cdots + P(W(r))C.
 \end{align*}
 Each term in this expansion contains matrix entries in the factors that share a common quaternion subalgebra, so we can use the simplification of $P(x)y =xyx$ available in associative cases for each term:
 \begin{align*}
  P(W(r))[x] = W(r) a W(r) + \cdots + W(r)C W(r) = W(r)[x]W(r).
 \end{align*}
 This confirms that for $r$ a commutative vector, the map $[x] \mapsto W(r) [x] W(r)$ is an isometry of $\mathbb{OP}^2$.
\end{proof}

\section{The Common Construction}

This section defines a common construction for pairs of projective $t$-designs, provides a familiar example, and then applies the common construction to the two strictly projective tight $5$-designs of Theorem~\ref{5designclassification}. 

\begin{definition}[Common construction]
 \label{commonconstruction}
 Let $r_1, r_2, \ldots, r_n$ be commutative row vectors in $\mathbb{F}^\rho$, with $d = [\mathbb{F}: \mathbb{R}]$, and let $[r_1], [r_2], \ldots, [r_n]$ be the corresponding primitive idempotents in projective space $\mathbb{FP}^{\rho-1}$.
 Let $G$ be the group acting on $\mathbb{F}^\rho$ generated by the following isometries under composition:
 \begin{align*}
  x \mapsto x\left(I_\rho - 2 [r_i]\right),~i = 1,2,\ldots, n.
 \end{align*}
 Let $H$ be the group acting on $\mathbb{FP}^{\rho- 1}$ generated by the following isometries under composition:
 \begin{align*}
  [x] \mapsto \left(I_\rho - 2 [r_i]\right)[x]\left(I_\rho - 2 [r_i]\right),~i = 1,2,\ldots, n.
 \end{align*}
 If $G$ is finite, then the orbit of $r_1,r_2 \ldots, r_n$ defines a spherical design in $\Omega_{\rho d}$ using $\mathbb{R}^{ \rho d} \cong \mathbb{F}^\rho$. 
 The lines spanned by the points of this spherical design define a real projective design in $\mathbb{R}^{\rho d-1}$.
 When $H$ is finite, the orbit of $[r_1], [r_2], \ldots, [r_n]$ defines a projective design in $\mathbb{FP}^{\rho-1}$.
\end{definition}

\begin{remark}
 This common construction definition relies on Lemma~\ref{isometrylemma} to ensure that the needed maps are indeed isometries when the vectors $r_1, r_2, \ldots, r_n$ are each commutative. 
 When $\mathbb{F} \ne \mathbb{O}$, we can relax the commutative requirement of Definition~\ref{commonconstruction}.
 In general, a selection $r_1, r_2, \ldots, r_n$ of commutative vectors must be chosen carefully in order for the isometries given in Definition~\ref{commonconstruction} to generate finite groups. 
 More details about the classification of finite reflection groups are available in~\cite{shephard_finite_1954, cohen_finite_1976, cohen_finite_1980, wilson_finite_2009}.
\end{remark}

\begin{example}
 \label{c6example}
 Let $r_1,r_2,\ldots, r_6 \in \mathbb{C}^6$ be the rows of the following matrix, where $\omega$ is a complex cube root of unity:
 \begin{align*}
  \begin{bmatrix}
    2 & 0 & 0 & 0 & 0 & 0 \\
    0 & 2 & 0 & 0 & 0 & 0 \\ 
    0 & 0 & 2 & 0 & 0 & 0 \\
    1 & \omega & \omega & 1 & 0 & 0 \\
    \omega & 1 & \omega & 0 & 1 & 0 \\
    \omega & \omega & 1 & 0 & 0 & 1 
  \end{bmatrix}.
 \end{align*}
 The common construction of Definition~\ref{commonconstruction} yields finite groups $G$ and $H$.
 The group $G$ acting on $\mathbb{C}^6$ is the complex reflection group $W(K_6) = (6.\PSU _4(3)) : 2$. The orbit of $r_1, r_2,\ldots, r_6$ under the action of $G$ form the $756$ shortest vectors of the $K_{12}$ integral lattice, defining a $5$-design in $\Omega_{12}$.
 This corresponds to the projective $2$-design in $\mathbb{RP}^{11}$ consisting of the $378$ lines spanned by the $K_{12}$ short vectors. 
 The orbit of $[r_1], [r_2], \ldots, [r_6]$ under the action of $H = \PSU _4(3) : 2$ form the $126$ points of a tight $3$-design in $\mathbb{CP}^5$. More details about this example are available in~\cite[127-129]{conway_sphere_2013}.
\end{example}

We may now introduce the main result of this paper.

\begin{theorem}
 There exists a set of commutative vectors $r_1, r_2, \ldots, r_n$ in $\mathbb{O}^3$ that yield the two strictly projective tight $5$-designs under the common construction of Definition~\ref{commonconstruction}.
 The vectors $r_1, r_2, \ldots, r_n$ are not unique and can be given, for example, by the rows of the following matrix, for any $t \in \mathbb{F}_7$ and with $s = \frac{1}{2}(-1+i_0 + i_1 + i_2 + i_3 + i_4 + i_5 + i_6)$:
 \begin{align*}
\begin{bmatrix}
    2 & 2 & 0 \\
    2s & 0 & 0 \\
    s^2 & s & s \\
    2 & 2i_{t} & 0 \\
    2 & 2i_{t+1} & 0 \\
    2 & 2i_{t+3} & 0 
\end{bmatrix}.
 \end{align*}
 \label{maintheorem}
\end{theorem}

\begin{proof}
 The examples given in Theorem~\ref{maintheorem} were found and checked using the software GAP~\cite{noauthor_gap_2020}. 
 Setting $t=0$, the computation begins by first applying the isometries of Definition~\ref{commonconstruction} respectively to $x \in \{r_1, \ldots, r_6\}$ and $[x] \in \{[r_1], \ldots, [r_6]\}$. Any new elements in $\mathbb{O}^3$ and $\mathbb{OP}^2$ are added to the respective sets. 
 This process is repeated and each application of the isometries either provides new elements or permutes the elements of the set.
 Once all six isometries permute the appropriate set, without providing new elements, those permutations are used to generate the groups $G$ and $H$. Permutation group tools in GAP identify $G$ and $H$ as $2 \cdot G_2(4)$ and $^3 D_4(2)$ respectively. The properties of the tight $5$-designs given by the full orbits of $\{r_1, \ldots, r_6\}$ and $\{[r_1], \ldots, [r_6]\}$ are verified directly using the definitions above. 
 Since $i_t \mapsto i_{t+1}$ is an automorphism of $\mathbb{O}$, the result is also true for $0 \ne t \in \mathbb{F}_7$. 
\end{proof}

\begin{remark}
 If we apply the common construction of Definition~\ref{commonconstruction} to the vectors given in Theorem~\ref{maintheorem} then we obtain an isometry group $G = 2\cdot G_2(4) \subset O(24)$ acting on $\mathbb{R}^{24} \cong \mathbb{O}^3$. 
 The orbit of $r_1, \ldots, r_6$ under the action of $G$ forms a $\sqrt{2} : 1$ scale copy of the short vectors of the Leech lattice, which define the unique tight $11$-design on $\Omega_{24}$ with cardinality $196560$. 
 The lines spanned by the spherical tight $11$-design vectors form the corresponding tight $5$-design in $\mathbb{RP}^{23}$ with cardinality $98280$. 
 The common construction also yields isometry group $^3 D_4(2) \subset F_4$ acting on $\mathbb{OP}^2$. 
 The orbit of $[r_1], [r_2], \ldots, [r_6]$ under the action of $H$ form a copy of the unique tight $5$-design on $\mathbb{OP}^2$ of cardinality $819$.
\end{remark}

\section{Leech Lattice Symmetries and the Octonion Projective Plane}

In light of Theorem~\ref{maintheorem}, we can use the involutionary isometries from our common construction to generate certain symmetries of the Leech lattice and exhibit their relation to the octonion projective plane. 
This section outlines a construction of the Suzuki chain of Leech lattice symmetries acting on $\mathbb{O}^3$ and describes their corresponding action on $\mathbb{OP}^2$ where possible. 

The vectors in $\mathbb{O}^3 \cong \mathbb{R}^{24}$ of the spherical tight $11$-design in Theorem~\ref{maintheorem} are precisely the short vectors of Wilson's octonion Leech lattice construction~\cite{wilson_octonions_2009, wilson_finite_2009, wilson_conways_2011}. 
In fact, a GAP computation confirms that any choice of $t \in \mathbb{F}_7$ in Theorem~\ref{maintheorem} yields the same orbit in $\mathbb{O}^3$ but distinct orbits in $\mathbb{OP}^2$ and distinct groups of type $2\cdot G_2(4)$ acting on $\mathbb{O}^3$. 
The union of these seven distinct $2\cdot G_2(4)$ groups generates the full Leech lattice automorphism group, Conway's group $\Co _0 = 2\cdot \Co _1$. 
The corresponding permutation group acting on the $98280$ lines defining the tight $5$-design in $\mathbb{RP}^{23}$ is the sporadic simple group $\Co _1$.

The group $\Co _1$ contains the alternating group $A_9$ as a subgroup, with symmetric group $S_3$ centralizing $A_9$ in $\Co _1$. The \textit{Suzuki chain} is a chain of centralizers in $\Co _1$ of the corresponding chain of alternating groups $A_9 > A_8 > \cdots > A_3$:
\begin{align*}
 S_3 < S_4 < \PSL _2(7) < \PSU _3(3) < \HJ < G_2(4) < 3\cdot \Suz .
\end{align*}
Here $\HJ $ and $\Suz $ are respectively the Hall--Janko and Suzuki sporadic simple groups. 
Wilson uses scalar octonion multiplication acting on a Leech lattice in $\mathbb{O}^3$ to construct the isometries needed to generate a chain of double covers of the alternating groups using in the Suzuki chain. 
Wilson then also includes known coordinate symmetries of the octonion triples (coordinate permutations and sign changes) to construct a maximal subgroup of $2\cdot \Co _1$. 
By appending the single reflection $x \mapsto x W(r)$, with $r = (s,1,1)$, Wilson is able to recover the entire Leech lattice automorphism group $\Co _0$. 
The Suzuki chain subgroups are described in the context of these isometries. More details about Wilson's construction and the group theory involved are available in~\cite{wilson_octonions_2009, wilson_finite_2009, wilson_conways_2011}.

Using the computational results of the construction in Theorem~\ref{maintheorem}, we can provide an alternative construction of Leech lattice automorphisms and the Suzuki chain. 
The benefit of this alternative construction is that all the generators involved are involutionary isometries of $\mathbb{O}^3$ with corresponding isometries on $\mathbb{OP}^2$ via the common construction of Definition~\ref{commonconstruction}. 

\begin{definition}
 Let $t \in \mathbb{F}_7$, let $s = \frac{1}{2}(-1 + i_0 + i_1 + i_2 + i_3 + i_4 + i_5 + i_6)$, and let $V_\infty$, $V_t$, and $S$ be, respectively, the sets of vectors of the form $(2,2,0)$, $(2,2 i_t, 0)$, and $(s^2, s, s)$, under all coordinate permutations and sign changes. 
\end{definition}

\begin{example}
\label{additionalexamples}
Under the common construction of Definition~\ref{commonconstruction}, the following commutative vectors yield the following finite group actions $G/\{\pm 1\}$ acting on $\mathbb{O}^3/\mathbb{R} \cong \mathbb{RP}^{23}$ and $H$ acting on $\mathbb{OP}^2$. Here we have $t,t'\in\mathbb{F}_7$ and $t\ne t'$.
\[
\belowdisplayskip12pt
\renewcommand{\arraycolsep}{8pt}
\renewcommand{\arraystretch}{1.2}
 \begin{array}{|c|c|c|}
\hline
 \{r_1, \ldots, r_n\}  & G/\{\pm 1\} \subset \Co _1 & H \subset F_4 \\
 \hline \hline 
 V_\infty & S_4 & S_4
 \\
 S & \PSL _2(7) & \PSL _2(7)
 \\
 S \cup V_t & \PSU _3(3) & \PSU _3(3)
 \\
 S \cup V_t \cup V_{t'} & \HJ & ^3 D_4(2)
 \\
 S \cup V_t \cup V_{t+1} \cup V_{t+3} & G_2(4) & ^3 D_4(2)
 \\
 S \cup V_{t+2} \cup V_{t+5} \cup V_{t+6} \cup V_{t+7} & 3\cdot \Suz & 
 \\
 S \cup V_{\infty} \cup V_{0} \cup V_{1} \cup V_{2} \cup V_{3} \cup V_{4} \cup V_{5} \cup V_{6} & \Co _1 & \\
 \hline
 \end{array}
\] 
\end{example}

\begin{remark}
 The groups described in Example~\ref{additionalexamples} are computed using GAP in the same manner as described in the proof of Theorem~\ref{maintheorem}.
 Computation time is saved by using the octonion automorphism $i_t \mapsto i_{t+1}$ to reduce the number of cases to check.
\end{remark}

\begin{remark}
 In Example~\ref{additionalexamples}, the bottom three rows yield the full tight $5$-design in $\mathbb{RP}^{23}$. Both the fourth and fifth rows yield the tight $5$-design in $\mathbb{OP}^2$. 
 The fifth row yields both tight $5$-designs and corresponds to the example in Theorem~\ref{maintheorem}.
\end{remark}

\begin{remark}
\looseness-1
 The first three rows of the table in Example~\ref{additionalexamples} involve initial vectors $\{r_1, \ldots, r_n\}$ with coefficients belonging to a common associative subalgebra of $\mathbb{O}$. Accordingly, the isometries generating $G$ and $H$ also generate matrix groups, which accounts for the agreement between $G/\{\pm 1\}$ and $H$. 
 The bottom four rows in the table involve sets $\{r_1, \ldots, r_n\}$ with coefficients that generate the full non-associative octonion algebra. This means that the groups $G$ and $H$ generated by the initial isometries are no longer related to a matrix group generated by matrices $W(r_i) = I_3 - 2[r_i]$. Indeed, the matrices $W(r_i)$ instead generate a non-associative octonion matrix loop rather than a matrix group. This partly explains the divergence between the properties of the groups $G$ and $H$ in the bottom four rows of the table.
 In the bottom two rows of the table, the initial isometries generate a finite group $G$ acting on $\mathbb{O}^3$ but do not seem to generate a corresponding finite group acting on $\mathbb{OP}^2$. 
 An open question is whether they generate the Lie group $F_4$, the full isometry group of $\mathbb{OP}^2$.
\end{remark}

\begin{remark}
\looseness-1
 As described above, Theorem~\ref{maintheorem} yields the same tight $5$-design in $\mathbb{RP}^{23}$ but distinct tight $5$-designs in $\mathbb{OP}^2$ for distinct $t \in \mathbb{F}_7$. 
 In contrast, if we take the conjugate $s \mapsto \overline{s}$ in Theorem~\ref{maintheorem}, then for distinct $t \in \mathbb{F}_7$ the common construction will instead yield distinct tight $5$-designs in $\mathbb{RP}^{23}$ but just one common tight $5$-design in $\mathbb{OP}^2$.
\end{remark}

\begin{remark}
 Possible variations on initial vectors in Theorem~\ref{maintheorem} include using carefully selected norm $2$ integral octonions to construct the initial vectors $\{r_1, \ldots, r_n\}$ of the common construction of Definition~\ref{commonconstruction}. A subsequent paper will explore constructions of this form and how they can be used to exhibit Suzuki chain symmetries of the Leech lattice. 
\end{remark}

\section{Conclusion}

We have seen in Theorem~\ref{maintheorem} and Example~\ref{additionalexamples} that the vectors $\{r_1, \ldots, r_n\} = S \cup V_{t} \cup V_{t+1} \cup V_{t+3}$, with $t \in \mathbb{F}_7$, define isometries of $\mathbb{O}^3$ and $\mathbb{OP}^2$ according to the common construction of Definition~\ref{commonconstruction}. 
These isometries generate a finite group $G_2(4)$ acting on $\mathbb{O}^3/\mathbb{R} \cong \mathbb{RP}^{23}$ and a finite group $^3 D_4(2)$ acting on $\mathbb{OP}^2$ that yield the two strictly projective tight $5$-designs as orbits. 
Specifically, the two tight $5$-designs are respectively the orbits of the initial vectors $\{r_1, \ldots, r_n\}$ and of the primitive idempotents $\{[r_1], \ldots, [r_n]\}$.
This common construction accounts for the previously conjectured connection between these two tight $5$-designs in~\cite{hoggar_t-designs_1982}.

Hoggar remarks in~\cite{hoggar_t-designs_1982} that his conjectured connection was met with skepticism. 
``Against this, the referee remarks: the automorphism group of the unique $(2,8)$ hexagon has index 3 subgroup $^3D_4(2)$, which has no irreducible projective representation of degree $\le 24$. Furthermore, $^3 D_4(2)$ has no proper subgroups acting transitively on the $819$ points''~\cite{hoggar_t-designs_1982}.
In our common construction of Definition~\ref{commonconstruction}, the link between our two designs is the initial commutative vectors $\{r_1, r_2, \ldots, r_n\}$ rather than the group $^3 D_4(2)$, its subgroups, or its representations. Indeed the pair of groups, $G_2(4)$ acting on $\mathbb{RP}^{23}$ and $^3D_4(2)$ acting on $\mathbb{OP}^2$, have relative cardinality $25/21$ so that one cannot be a subgroup of the other. The non-associativity of the octonion algebra permits the common construction to yield these seemingly unrelated groups and the tight 5-designs given by their orbits.

\subsubsection*{Acknowledgements} The author thanks his two dissertation supervisors, Claude Tardif and Charles Paquette, for their advice and input. The author also thanks an anonymous reviewer for their helpful suggestions. During this research the author received funding from an Individualized Learning Plan, through the Canadian Defence Academy, and also from the Royal Military College Bursary Fund.

% \nocite{*}

\bibliographystyle{amsalpha}
\bibliography{Bibliography}

\end{document}